\documentclass[twoside,reqno,11pt]{amsart}

\usepackage{latexsym}

\newcommand{\qdn}{\hspace*{-1.5mm}}
\newcommand{\qqdn}{\hspace*{-2.5mm}}
\newcommand{\xqdn}{\hspace*{-5.0mm}}
\newcommand{\xxqdn}{\hspace*{-10mm}}

\newcommand{\fns}{\footnotesize}

\newcommand{\ffnk}[4]{\left[\qdn\ba{#1}#3\\#4\ea{\!\Big|\:#2}\right]}

\newcommand{\binm}{\binom}

\newcommand{\binq}[2]{\genfrac{[}{]}{0mm}{0}{#1}{#2}}

\newcommand{\nnm}{\nonumber}
\newcommand{\be}{\begin{equation}}
\newcommand{\ee}{\end{equation}}
\newcommand{\ba}{\begin{array}}
\newcommand{\ea}{\end{array}}
\newcommand{\bmn}{\begin{eqnarray}}
\newcommand{\emn}{\end{eqnarray}}
\newcommand{\bnm}{\begin{eqnarray*}}
\newcommand{\enm}{\end{eqnarray*}}
\newcommand{\bln}{\begin{subequations}}
\newcommand{\eln}{\end{subequations}}

\newtheorem{thm}{Theorem}

\newtheorem{corl}[thm]{Corollary}

\newtheorem{entry}{Entry}

\newcommand{\bbtm}[4]{\bibitem{kn:#1}{#2,}~{#3,}~{#4.}}
\newcommand{\cito}[1]{\cite{kn:#1}}
\newcommand{\citu}[2]{\cite[#2]{kn:#1}}

\begin{document} 
{\fns
\title{Summation formulae for $q$-Watson\\ type $_4\phi_3$-series}

\author{$^a$Chuanan Wei, $^b$Dianxuan Gong, $^c$ Jianbo Li}
\dedicatory{
$^A$Department of Information Technology\\
  Hainan Medical College, Haikou 571101, China\\
   $^B$College of Sciences\\
   Hebei Polytechnic University, Tangshan 063009, China\\\
  $^C$Institute of Mathematical Sciences\\
 Xuzhou Normal University, Xuzhou 221116, China}
\thanks{}

\address{ }


\begin{abstract}
According to the method of series rearrangement, we establish two
families of summation formulae for $q$-Watson type
$_4\phi_3$-series.
\end{abstract}

\maketitle\thispagestyle{empty}
\markboth{Chuanan Wei, Dianxuan Gong, Jianbo Li}
         {Summation formulae for $q$-Watson type $_4\phi_3$-series}

\section{Introduction}

 For two complex numbers $x$ and
$q$, define the $q$-shifted factorial by
 \[(x;q)_0=1\quad\text{and}\quad (x;q)_n=\prod_{i=0}^{n-1}(1-xq^i)
  \quad\text{when}\quad n\in \mathbb{N}.\]
Its fraction form reads as
\[\qqdn\qdn\ffnk{ccccc}{q}{a,&b,&\cdots,&c}{\alpha,&\beta,&\cdots,&\gamma}_n
=\frac{(a;q)_n(b;q)_n\cdots(c;q)_n}{(\alpha;q)_n(\beta;q)_n\cdots(\gamma;q)_n}.\]
 Following Gasper and Rahman \cito{gasper}, the basic hypergeometric series
can be defined by
\[_{1+r}\phi_s\ffnk{cccc}{q;z}{a_0,&a_1,&\cdots,&a_r}
{&b_1,&\cdots,&b_s}
 =\sum_{k=0}^\infty
\ffnk{ccccc}{q}{a_0,&a_1,&\cdots,&a_r}{q,&b_1,&\cdots,&b_r}_kz^k\]
where $\{a_i\}_{i\geq0}$ and $\{b_j\}_{j\geq1}$ are complex
parameters such that no zero factors appear in the denominators of
the summand on the right hand side. Then the $q$-Watson formula due
to Andrews \cito{andrews} and the $q$-Watson formula due to Jain
\citu{jain}{Eq. 3.17} can be stated, respectively, as follows:
  \bmn\qquad \label{andrews}
{_4\phi_3}\ffnk{ccccc}{q;q}{q^{-n},q^{1+n}a,\sqrt{c},-\sqrt{c}}
{q\sqrt{a},-q\sqrt{a},c}=\begin{cases}
c^s\ffnk{ccccc}{q^2}{q,q^2a/c}{q^{2}a,qc}_s,&\qqdn n=2s;
\\[4mm]
0,&\qqdn n=1+2s,
\end{cases}
 \emn
  \bmn \label{jain}
{_4\phi_3}\ffnk{ccccc}{q;q}{a,c,q^{-n},-q^{-n}}
{\sqrt{qac},-\sqrt{qac},q^{-2n}}=\ffnk{ccccc}{q^2}{qa,qc}{q,qac}_n.
 \emn

 Inspired by the method due to Chu \cito{chu}, we shall establish
the following two families of summation formulae for $q$-Watson type
$_4\phi_3$-series:
\[{_4\phi_3}\ffnk{ccccc}{q;q}{q^{-n},q^{1+n}a,\sqrt{c},-\sqrt{c}}
{q\sqrt{a},-q\sqrt{a},q^\varepsilon c}\quad\text{and}\quad
{_4\phi_3}\ffnk{ccccc}{q;q}{q^{-n},q^{1+n}a,q^\varepsilon\!\sqrt{c},-\sqrt{c}}
{q\sqrt{a},-q\sqrt{a},q^\varepsilon c},\]
\[\xqdn{_4\phi_3}\ffnk{ccccc}{q;q}{a,c,q^{-n},-q^{-n}}
{\sqrt{qac},-\sqrt{qac},q^{\varepsilon-2n}}\quad\:\,\text{and}\quad
{_4\phi_3}\ffnk{ccccc}{q;q}{a,c,q^{\varepsilon-n},-q^{-n}}
{\sqrt{qac},-\sqrt{qac},q^{\varepsilon-2n}}\]
 where the disturbing parameter $\varepsilon$ is a nonnegative integer throughout the paper.
Eight concrete summation formulae corresponding
$1\leq\varepsilon\leq2$ will be exhibited.
\section{Summation formulae for $q$-Andrews-Watson type $_4\phi_3$-series}
 \subsection{} Letting $a\to c/q$, $b\to q^{-k}$ and $c\to \infty$ for the
terminating $_6\phi_5$-series identity (cf. Gasper and Rahman
\citu{gasper}{p. 42}):
 \bmn\label{terminating-65}
\qquad {_6\phi_5}\ffnk{cccccccc}{q;\frac{q^{1+\varepsilon}a}{bc}}
{a,\:q\sqrt{a},\:-q\sqrt{a},\:b,\:c,\:q^{-\varepsilon}}
 {\sqrt{a},\:-\sqrt{a},\:qa/b,\:qa/c,\:q^{1+\varepsilon}a}
 =\ffnk{ccccc}{q}{qa,qa/bc}{qa/b,qa/c}_\varepsilon,
 \emn
we obtain the following equation:
\[\sum_{i=0}^\varepsilon\binq{k}{i}q^{(i+\varepsilon-1)i}c^{i}\!\!
\ffnk{ccccc}{q}{cq^{k+i}}{cq^i}_{\varepsilon-i}\qdn
\frac{1-cq^{2i-1}}{1-cq^{i-1}}
\ffnk{ccccc}{q}{q^{-\varepsilon}}{q^{\varepsilon}c}_i=1.\]
 Then we can proceed as follows:
  \bnm
&&\xqdn{_4\phi_3}\ffnk{ccccc}{q;q}{q^{-n},q^{1+n}a,\sqrt{c},-\sqrt{c}}
{q\sqrt{a},-q\sqrt{a},q^{\varepsilon}c}=\sum_{k=0}^n(-1)^k\binq{n}{k}q^{\binm{k+1}{2}-nk}
\ffnk{ccccc}{q}{q^{1+n}a,\sqrt{c},-\sqrt{c}}{q\sqrt{a},-q\sqrt{a},q^{\varepsilon}c}_k\\
&&\xqdn\:\:=\:\sum_{k=0}^n(-1)^k\binq{n}{k}q^{\binm{k+1}{2}-nk}
\ffnk{ccccc}{q}{q^{1+n}a,\sqrt{c},-\sqrt{c}}{q\sqrt{a},-q\sqrt{a},q^{\varepsilon}c}_k\\
&&\xqdn\:\:\times\,\,
\sum_{i=0}^\varepsilon\binq{k}{i}q^{(i+\varepsilon-1)i}c^{i}\!\!
\ffnk{ccccc}{q}{cq^{k+i}}{cq^i}_{\varepsilon-i}\qdn
\frac{1-cq^{2i-1}}{1-cq^{i-1}}
\ffnk{ccccc}{q}{q^{-\varepsilon}}{q^{\varepsilon}c}_i.
 \enm
Interchanging the summation order, we can manipulate the last double
sum as
 \bnm
{_4\phi_3}\ffnk{ccccc}{q;q}{q^{-n},q^{1+n}a,\sqrt{c},-\sqrt{c}}
{q\sqrt{a},-q\sqrt{a},q^{\varepsilon}c}
&&\xqdn=\sum_{i=0}^\varepsilon\binq{n}{i}q^{(i+\varepsilon-1)i}c^i
\frac{1-cq^{2i-1}}{1-cq^{i-1}}
\ffnk{ccccc}{q}{q^{-\varepsilon}}{q^{\varepsilon}c}_i
\\&&\xqdn\times\:
\sum_{k=i}^n(-1)^k\binq{n-i}{k-i}q^{\binm{k+1}{2}-nk}
\ffnk{ccccc}{q}{q^{1+n}a,\sqrt{c},-\sqrt{c}}{q\sqrt{a},-q\sqrt{a},q^ic}_k.
 \enm
Shifting the summation index $k\to i+j$ for the sum on the last
line, we get the relation:
 \bmn
 {_4\phi_3}\ffnk{ccccc}{q;q}{q^{-n},q^{1+n}a,\sqrt{c},-\sqrt{c}}
{q\sqrt{a},-q\sqrt{a},q^{\varepsilon}c}&&\xqdn=\sum_{i=0}^{\varepsilon}
q^{(i+\varepsilon)i}c^i
\ffnk{ccccc}{q}{q^{-\varepsilon},q^{-n},q^{1+n}a,\sqrt{c},-\sqrt{c}}
{q,q\sqrt{a},-q\sqrt{a},q^\varepsilon c,q^{i-1}c}_i
        \nnm  \\&&\xqdn\times\:
{_4\phi_3}\ffnk{ccccc}{q;q}{q^{i-n},q^{1+n+i}a,q^i\sqrt{c},-q^i\sqrt{c}}
{q^{1+i}\!\sqrt{a},-q^{1+i}\!\sqrt{a},q^{2i}c}.\label{equation-a}
 \emn
 Evaluating the
$_4\phi_3$-series on the last line by \eqref{andrews}, we establish
the following summation theorem.

\begin{thm}\label{thm-a}
For two complex numbers $\{a,c\}$ and a nonnegative integer
$\varepsilon$, there holds the summation formula for
$q$-Andrews-Watson type $_4\phi_3$-series:
 \bnm
{_4\phi_3}\ffnk{ccccc}{q;q}{q^{-n},q^{1+n}a,\sqrt{c},-\sqrt{c}}
{q\sqrt{a},-q\sqrt{a},q^{\varepsilon}c}&&\xqdn=\sum_{i=0}^{\varepsilon}
q^{(\varepsilon+n)i}c^{\frac{n+i}{2}}
\ffnk{ccccc}{q}{q^{-\varepsilon},q^{-n},q^{1+n}a,\sqrt{c},-\sqrt{c}}
{q,q\sqrt{a},-q\sqrt{a},q^\varepsilon c,q^{i-1}c}_i
 \\&&\xqdn\times\:
\ffnk{ccccc}{q^2}{q,q^2a/c}{q^{2+2i}a,q^{1+2i}c}_{\frac{n-i}{2}}\chi(n-i\equiv_20)
 \enm
where $n-i\equiv_20$ stands for the congruence relation $n-i=0$ (mod
2) and $\chi$ is the logical function defined by
$\chi(\text{true})=1$ and $\chi(\text{false})=0$ otherwise.
\end{thm}

\begin{corl}[$\varepsilon=1$ in Theorem \ref{thm-a}]
\bnm
 \quad{_4\phi_3}\ffnk{ccccc}{q;q}{q^{-n},q^{1+n}a,\sqrt{c},-\sqrt{c}}
{q\sqrt{a},-q\sqrt{a},qc}\!=\!
\begin{cases}
c^s\ffnk{ccccc}{q^2}{q,q^2a/c}{q^{2}a,qc}_s,&\qqdn n=2s;
\\[4mm]
c^{1+s}\ffnk{ccccc}{q^2}{q}{qc}_{1+s}\!\!
 \ffnk{ccccc}{q^2}{q^2a/c}{q^2a}_{s},&\qqdn n=1+2s.
\end{cases}
\enm
\end{corl}

\begin{corl}[$\varepsilon=2$ in Theorem \ref{thm-a}]
\bnm
 &&\xxqdn\qdn{_4\phi_3}\ffnk{ccccc}{q;q}{q^{-n},q^{1+n}a,\sqrt{c},-\sqrt{c}}
{q\sqrt{a},-q\sqrt{a},q^2c}\\[2mm]&&\xxqdn=
\begin{cases}
\Big\{1+\tfrac{q^2c(1-c)(1-q^{2s})(1-q^{1+2s}a)}{(1-q^2c)(1-q^{1+2s}c)(1-q^{2s}a/c)}\Big\}
c^s\ffnk{ccccc}{q^2}{q,q^2a/c}{q^{2}a,qc}_s,&\qqdn n=2s;
\\[4mm]
\tfrac{1-q^2}{1-q^2c}\,c^{1+s}
\ffnk{ccccc}{q^2}{q^3,q^2a/c}{q^2a,q^3c}_s,&\qqdn n=1+2s.
\end{cases}
\enm
\end{corl}
\subsection{}

Letting $a\to c/q$, $b\to q^{-k}$ and $c\to \sqrt{c}$ for
\eqref{terminating-65}, we attain the following equation:
 \bnm
&&\sum_{i=0}^\varepsilon(-1)^iq^{\varepsilon+\binm{i}{2}}c^{\frac{i-\varepsilon}{2}}\frac{1-cq^{2i-1}}{1-c\,q^{i+\varepsilon-1}}
\ffnk{ccccc}{q}{q^{-\varepsilon}}{q}_{i}\ffnk{ccccc}{q}{q^{1-\varepsilon}/\sqrt{c}}{q^{2-\varepsilon-i}/c}_{\varepsilon}
\\&&\:\,\times\:\frac{<q^k;q>_i\:<cq^{k+\varepsilon-1};q>_{\varepsilon-i}}{(q^k\sqrt{c};q)_\varepsilon}
=1
 \enm
where we have used the symbol:
 \[<x;q>_0=1\quad\text{and}\quad <x;q>_n=\prod_{i=0}^{n-1}(1-xq^{-i})
  \quad\text{when}\quad n\in \mathbb{N}.\]
 Then we can proceed as follows:
  \bnm
&&\xqdn{_4\phi_3}\ffnk{ccccc}{q;q}{q^{-n},q^{1+n}a,q^{\varepsilon}\!\sqrt{c},-\sqrt{c}}
{q\sqrt{a},-q\sqrt{a},q^{\varepsilon}c}=\sum_{k=0}^n\ffnk{ccccc}{q;q}{q^{-n},q^{1+n}a,q^{\varepsilon}\!\sqrt{c},-\sqrt{c}}
{q,q\sqrt{a},-q\sqrt{a},q^{\varepsilon}c}_kq^k\\
&&\xqdn\:\:=\:\sum_{k=0}^n\ffnk{ccccc}{q;q}{q^{-n},q^{1+n}a,q^{\varepsilon}\!\sqrt{c},-\sqrt{c}}
{q,q\sqrt{a},-q\sqrt{a},q^{\varepsilon}c}_kq^k
\sum_{i=0}^\varepsilon(-1)^iq^{\varepsilon+\binm{i}{2}}c^{\frac{i-\varepsilon}{2}}\frac{1-cq^{2i-1}}{1-cq^{i+\varepsilon-1}}
\\&&\xqdn\:\:\times\:\ffnk{ccccc}{q}{q^{-\varepsilon}}{q}_{i}\ffnk{ccccc}{q}{q^{1-\varepsilon}/\sqrt{c}}{q^{2-\varepsilon-i}/c}_{\varepsilon}
\frac{<q^k;q>_i\:<cq^{k+\varepsilon-1};q>_{\varepsilon-i}}{(q^k\sqrt{c};q)_\varepsilon}.
 \enm
Interchanging the summation order, we can reformulate the last
double sum as
  \bnm
&&\xqdn{_4\phi_3}\ffnk{ccccc}{q;q}{q^{-n},q^{1+n}a,q^{\varepsilon}\!\sqrt{c},-\sqrt{c}}
{q\sqrt{a},-q\sqrt{a},q^{\varepsilon}c}=
\sum_{i=0}^\varepsilon(-1)^iq^{\varepsilon+\binm{i}{2}}c^{\frac{i-\varepsilon}{2}}\frac{1-cq^{2i-1}}{1-cq^{i+\varepsilon-1}}
\\&&\xqdn\:\:\times\:\ffnk{ccccc}{q}{q^{-\varepsilon}}{q}_{i}\ffnk{ccccc}{q}{q^{1-\varepsilon}/\sqrt{c}}{q^{2-\varepsilon-i}/c}_{\varepsilon}
\\&&\xqdn\:\:\times\: \sum_{k=i}^n\ffnk{ccccc}{q;q}{q^{-n},q^{1+n}a,q^{\varepsilon}\!\sqrt{c},-\sqrt{c}}
{q,q\sqrt{a},-q\sqrt{a},q^{\varepsilon}c}_kq^k
\frac{<q^k;q>_i\:<cq^{k+\varepsilon-1};q>_{\varepsilon-i}}{(q^k\sqrt{c};q)_\varepsilon}.
 \enm
Shifting the summation index $k\to i+j$ for the sum on the last
line, we achieve the relation:
  \bmn
&&{_4\phi_3}\ffnk{ccccc}{q;q}{q^{-n},q^{1+n}a,q^{\varepsilon}\!\sqrt{c},-\sqrt{c}}
{q\sqrt{a},-q\sqrt{a},q^{\varepsilon}c}\nnm\\&&\:\:=\!\:
\sum_{i=0}^\varepsilon(-1)^i
 q^{\varepsilon i+\binm{i+1}{2}}c^{\frac{i}{2}}\ffnk{ccccc}{q}{q^{-\varepsilon},q^{-n},q^{1+n}a,\sqrt{c},-\sqrt{c}}
{q,q\sqrt{a},-q\sqrt{a},q^\varepsilon c,q^{i-1}c}_i
\nnm\\&&\:\:\times\,\:{_4\phi_3}\ffnk{ccccc}{q;q}{q^{i-n},q^{1+n+i}a,q^i\!\sqrt{c},-q^i\sqrt{c}}
{q^{1+i}\sqrt{a},-q^{1+i}\sqrt{a},q^{2i}c}. \label{equation-b}
 \emn
 Evaluating the
$_4\phi_3$-series on the last line by \eqref{andrews}, we found the
following summation theorem.

\begin{thm}\label{thm-b}
For two complex numbers $\{a,c\}$ and a nonnegative integer
$\varepsilon$, there holds the summation formula for
$q$-Andrews-Watson type $_4\phi_3$-series:
 \bnm
&&\xqdn{_4\phi_3}\ffnk{ccccc}{q;q}{q^{-n},q^{1+n}a,q^{\varepsilon}\!\sqrt{c},-\sqrt{c}}
{q\sqrt{a},-q\sqrt{a},q^{\varepsilon}c}\\&&\xqdn\:\:=\:\sum_{i=0}^{\varepsilon}
(-1)^iq^{(\varepsilon+n)i-\binm{i}{2}}c^{\frac{n}{2}}
\ffnk{ccccc}{q}{q^{-\varepsilon},q^{-n},q^{1+n}a,\sqrt{c},-\sqrt{c}}
{q,q\sqrt{a},-q\sqrt{a},q^\varepsilon c,q^{i-1}c}_i
 \\&&\xqdn\:\:\times\:
\ffnk{ccccc}{q^2}{q,q^2a/c}{q^{2+2i}a,q^{1+2i}c}_{\frac{n-i}{2}}\chi(n-i\equiv_20).
 \enm
\end{thm}

\begin{corl}[$\varepsilon=1$ in Theorem \ref{thm-b}]
\bnm
 \quad{_4\phi_3}\ffnk{ccccc}{q;q}{q^{-n},q^{1+n}a,q\sqrt{c},-\sqrt{c}}
{q\sqrt{a},-q\sqrt{a},qc}\!=\!
\begin{cases}
c^s\ffnk{ccccc}{q^2}{q,q^2a/c}{q^{2}a,qc}_s,&\qqdn n=2s;
\\[4mm]
-c^{\frac{1}{2}+s}\!\ffnk{ccccc}{q^2}{q}{qc}_{1+s}\!\!
 \ffnk{ccccc}{q^2}{q^2a/c}{q^2a}_{s}\!,&\qqdn n=1+2s.
\end{cases}
\enm
\end{corl}

\begin{corl}[$\varepsilon=2$ in Theorem \ref{thm-b}]
\bnm
 &&\xqdn\qdn{_4\phi_3}\ffnk{ccccc}{q;q}{q^{-n},q^{1+n}a,q^2\!\sqrt{c},-\sqrt{c}}
{q\sqrt{a},-q\sqrt{a},q^2c}\\[2mm]&&\xqdn=
\begin{cases}
\Big\{1+\tfrac{q(1-c)(1-q^{2s})(1-q^{1+2s}a)}{(1-q^2c)(1-q^{1+2s}c)(1-q^{2s}a/c)}\Big\}
\,c^s\ffnk{ccccc}{q^2}{q,q^2a/c}{q^{2}a,qc}_s,&\qqdn n=2s;
\\[4mm]
\tfrac{q^2-1}{1-q^2c}\,c^{\frac{1}{2}+s}
\ffnk{ccccc}{q^2}{q^3,q^2a/c}{q^2a,q^3c}_s,&\qqdn n=1+2s.
\end{cases}
\enm
\end{corl}
\section{Summation formulae for $q$-Jain-Watson type $_4\phi_3$-series}
\subsection{}

Performing the replacement $a\to aq^{-1-n}$ for \eqref{equation-a},
we have
 \bnm
 \,{_4\phi_3}\ffnk{ccccc}{\!q;q}{q^{-n},a,\sqrt{c},-\sqrt{c}}
{\sqrt{q^{1-n}a},-\sqrt{q^{1-n}a},q^{\varepsilon}c}&&\!\xqdn\!=\!\sum_{i=0}^{\varepsilon}\!
q^{(i+\varepsilon)i}c^i\!
\ffnk{ccccc}{q\!}{q^{-\varepsilon},q^{-n},a,\sqrt{c},-\sqrt{c}}
{q,\sqrt{q^{1-n}a},-\sqrt{q^{1-n}a},q^\varepsilon c,q^{i-1}c}_i
        \nnm  \\&&\!\xqdn\times\:
{_4\phi_3}\ffnk{ccccc}{q;q}{q^{i-n},q^{i}a,q^i\!\sqrt{c},-q^i\!\sqrt{c}}
{q^{i}\!\sqrt{q^{1-n}a},-q^{i}\!\sqrt{q^{1-n}a},q^{2i}c}.
 \enm
Employing the substitutions $c\to q^{-2n}$ and $q^{-n}\to c$ for the
last equation, we obtain the following relation:
 \bnm
 {_4\phi_3}\ffnk{ccccc}{q;q}{a,c,q^{-n},-q^{-n}}
{\sqrt{qac},-\sqrt{qac},q^{\varepsilon-2n}}&&\!\!\xqdn=\!\sum_{i=0}^{\varepsilon}
q^{(i+\varepsilon-2n)i}
\ffnk{ccccc}{q}{q^{-\varepsilon},q^{-n},-q^{-n},a,c}
{q,\sqrt{qac},-\sqrt{qac},q^{\varepsilon-2n},q^{i-1-2n}}_i
        \nnm  \\&&\!\xqdn\times\:
{_4\phi_3}\ffnk{ccccc}{q;q}{q^{i}a,q^{i}c,q^{i-n},-q^{i-n}}
{q^{i}\!\sqrt{qac},-q^{i}\!\sqrt{qac},q^{2i-2n}}.
 \enm
 Evaluating the
$_4\phi_3$-series on the last line by \eqref{jain}, we establish the
following summation theorem.

\begin{thm}\label{thm-c}
For two complex numbers $\{a,c\}$ and a nonnegative integer
$\varepsilon$ with $0\leq\varepsilon\leq n$, there holds the
summation formula for $q$-Jain-Watson type $_4\phi_3$-series:
 \bnm
\:{_4\phi_3}\ffnk{ccccc}{q;q}{a,c,q^{-n},-q^{-n}}
{\sqrt{qac},-\sqrt{qac},q^{\varepsilon-2n}}&&\!\xqdn=\!\sum_{i=0}^{\varepsilon}
q^{(i+\varepsilon-2n)i}
\ffnk{ccccc}{q}{q^{-\varepsilon},q^{-n},-q^{-n},a,c}
{q,\sqrt{qac},-\sqrt{qac},q^{\varepsilon-2n},q^{i-1-2n}}_i
        \nnm  \\&&\xqdn\!\times\:
\ffnk{ccccc}{q^2}{q^{1+i}a,q^{1+i}c}{q,q^{1+2i}ac}_{n-i}.
 \enm
\end{thm}

\begin{corl}[$\varepsilon=1$ in Theorem \ref{thm-c}: $n\geq1$]
\bnm \:\:\: {_4\phi_3}\ffnk{ccccc}{q;q}{a,c,q^{-n},-q^{-n}}
{\sqrt{qac},-\sqrt{qac},q^{1-2n}}=
\ffnk{ccccc}{q^2}{qa,qc}{q,qac}_{n}+\ffnk{ccccc}{q^2}{a,c}{q,qac}_{n}.
 \enm
\end{corl}

\begin{corl}[$\varepsilon=2$ in Theorem \ref{thm-c}: $n\geq2$]
\bnm
 {_4\phi_3}\ffnk{ccccc}{q;q}{a,c,q^{-n},-q^{-n}}
{\sqrt{qac},-\sqrt{qac},q^{2-2n}}&&\xqdn\!=
\frac{(1+q)(1-q^{1-2n})}{1-q^{2-2n}}\ffnk{ccccc}{q^2}{a,c}{q,qac}_{n}\\
&&\xxqdn\xxqdn\xqdn\qqdn+\Big\{1\!+\!\tfrac{q(1-a)(1-c)(1-q^{-2n})}{(1-q^{2-2n})(1-aq^{2n-1})(1-cq^{2n-1})}\Big\}\!
\ffnk{ccccc}{q^2}{qa,qc}{q,qac}_{n}.
 \enm
\end{corl}

\subsection{}
Performing the replacement $a\to aq^{-1-n}$ for \eqref{equation-b},
we have
 \bnm
  &&{_4\phi_3}\ffnk{ccccc}{\!q;q}{q^{-n},a,q^{\varepsilon}\!\sqrt{c},-\sqrt{c}}
{\sqrt{q^{1-n}a},-\sqrt{q^{1-n}a},q^{\varepsilon}c}\\&&\:\:=\:\sum_{i=0}^{\varepsilon}
(-1)^iq^{\varepsilon i+\binm{i+1}{2}}c^{\frac{i}{2}}\!
\ffnk{ccccc}{q\!}{q^{-\varepsilon},q^{-n},a,\sqrt{c},-\sqrt{c}}
{q,\sqrt{q^{1-n}a},-\sqrt{q^{1-n}a},q^\varepsilon c,q^{i-1}c}_i
        \nnm  \\&&\:\:\times\:\,
{_4\phi_3}\ffnk{ccccc}{q;q}{q^{i-n},q^{i}a,q^i\sqrt{c},-q^i\sqrt{c}}
{q^{i}\!\sqrt{q^{1-n}a},-q^{i}\!\sqrt{q^{1-n}a},q^{2i}c}.
 \enm
Employing the substitutions $c\to q^{-2n}$ and $q^{-n}\to c$ for the
last equation, we get the following relation:
 \bnm
 &&{_4\phi_3}\ffnk{ccccc}{q;q}{a,c,q^{\varepsilon-n},-q^{-n}}
{\sqrt{qac},-\sqrt{qac},q^{\varepsilon-2n}}\\&&\:\:=\:\sum_{i=0}^{\varepsilon}
(-1)^iq^{(\varepsilon-n)i+\binm{i+1}{2}}
\ffnk{ccccc}{q}{q^{-\varepsilon},q^{-n},-q^{-n},a,c}
{q,\sqrt{qac},-\sqrt{qac},q^{\varepsilon-2n},q^{i-1-2n}}_i
        \nnm  \\&&\:\:\times\:\,
{_4\phi_3}\ffnk{ccccc}{q;q}{q^{i}a,q^{i}c,q^{i-n},-q^{i-n}}
{q^{i}\!\sqrt{qac},-q^{i}\!\sqrt{qac},q^{2i-2n}}.
 \enm

Evaluating the $_4\phi_3$-series on the last line by \eqref{jain},
we found the following summation theorem.

\begin{thm}\label{thm-d}
For two complex numbers $\{a,c\}$ and a nonnegative integer
$\varepsilon$ with $0\leq\varepsilon\leq n$, there holds the
summation formula for $q$-Jain-Watson type $_4\phi_3$-series:
 \bnm
&&{_4\phi_3}\ffnk{ccccc}{q;q}{a,c,q^{\varepsilon-n},-q^{-n}}
{\sqrt{qac},-\sqrt{qac},q^{\varepsilon-2n}}\\&&\:\:=\:\sum_{i=0}^{\varepsilon}
(-1)^iq^{(\varepsilon-n)i+\binm{i+1}{2}}
\ffnk{ccccc}{q}{q^{-\varepsilon},q^{-n},-q^{-n},a,c}
{q,\sqrt{qac},-\sqrt{qac},q^{\varepsilon-2n},q^{i-1-2n}}_i
        \nnm  \\&&\:\:\times\:
\ffnk{ccccc}{q^2}{q^{1+i}a,q^{1+i}c}{q,q^{1+2i}ac}_{n-i}.
 \enm
\end{thm}

\begin{corl}[$\varepsilon=1$ in Theorem \ref{thm-d}: $n\geq1$]
\bnm
 \qquad{_4\phi_3}\ffnk{ccccc}{q;q}{a,c,q^{1-n},-q^{-n}}
{\sqrt{qac},-\sqrt{qac},q^{1-2n}}=
\ffnk{ccccc}{q^2}{qa,qc}{q,qac}_{n}-q^n\ffnk{ccccc}{q^2}{a,c}{q,qac}_{n}.
 \enm
\end{corl}

\begin{corl}[$\varepsilon=2$ in Theorem \ref{thm-d}: $n\geq2$]
\bnm
 {_4\phi_3}\ffnk{ccccc}{q;q}{a,c,q^{-n},-q^{-n}}
{\sqrt{qac},-\sqrt{qac},q^{2-2n}}&&\xqdn\!=
\frac{(1+q)(1-q^{2n-1})}{q^{n-1}-q^{1-n}}\ffnk{ccccc}{q^2}{a,c}{q,qac}_{n}\\
&&\xxqdn\xxqdn\xqdn\qqdn\!+\Big\{1\!-\!\tfrac{(1-a)(1-c)(1-q^{2n})}{(1-q^{2-2n})(1-aq^{2n-1})(1-cq^{2n-1})}\Big\}\!
\ffnk{ccccc}{q^2}{qa,qc}{q,qac}_{n}.
 \enm
\end{corl}

With the change of $\varepsilon$, more concrete formulae could be
derived from Theorems \ref{thm-a}, \ref{thm-b}, \ref{thm-c} and
\ref{thm-d}. Considering that the resulting identities will become
complicated, we shall not display them one by one. Altough Chu
\cito{chu} gave more formulae for hypergeometric series,  we end our
paper in the present form because of the restriction of basic
hypergeometric series.


\end{document}